\documentstyle[12pt]{article}

\newtheorem{thm}{Theorem}
\newtheorem{lem}[thm]{Lemma}

\newtheorem{pr}[thm]{Proposition}

\newenvironment{mainlem}{\vspace{2ex}\noindent{\bf Main Lemma} \em }{\\}
\newenvironment{pf}{\noindent {\em Proof:}}{$\Box$\\}

\newcommand{\N}{\mbox{\hskip.1em N \hskip -1.25em \relax I \hskip .1em}}
\newcommand{\R}{\mbox{\hskip.1em R \hskip -1.25em \relax I \hskip .1em}}
\newcommand{\sgn}{\mbox{sgn\,}}
\newcommand{\spn}{\mbox{span}}
\newcommand{\sumi}{\sum^{i}_{j=1}}
\newcommand{\sumk}{\sum^k_{i=1}}
\newcommand{\sumsv}{\sum_{i>4}}
\newcommand{\rnorm}{\biggr\|}
\newcommand{\lnorm}{\biggl\|}

\newcommand{\pinf}{{p,\infty}}
\newcommand{\ep}{\epsilon^{1/4}}

\begin{document}

\begin{center}
  {\Large\bf  On $c_0$-saturated Banach spaces}\vspace{3mm}\\
  {\large\sc Denny H. Leung}
\end{center}

\vspace{1mm}

\begin{abstract}
  A Banach space $E$\/ is $c_0$-saturated if every closed infinite
dimensional subspace of $E$ contains an isomorph of $c_0$.  A
$c_0$-saturated Banach space with an unconditional basis which has a
quotient space isomorphic to $\ell^2$ is constructed.
\end{abstract}

\begin{figure}[b]
  \rule{3in}{.005in}\\
  1991 {\em Mathematics Subject Classification\ } 46B10.
\end{figure}

\baselineskip 4ex

 A Banach space $E$\/ is $c_0$-saturated if every closed infinite
dimensional subspace of $E$ contains an isomorph of $c_0$. In
\cite{O} and \cite{R}, it was asked whether all quotient spaces of
$c_0$-saturated spaces having  unconditional bases are also
$c_0$-saturated.  In \cite{R}, Rosenthal expressed the opinion that
the answer should be no.  Here, we construct an example which confirms
this opinion.\\
\indent Standard Banach space terminology, as may be found in
\cite{LT}, is employed.  For $1 \leq p \leq \infty$, $\|\cdot\|_p$
denotes the $\ell^p$-norm, and if additionally $p < \infty$,
\[ \|(a_n)\|_\pinf = \sup a^*_nn^{1/p}, \]
where $(a^*_n)$ is the decreasing rearrangement of $(|a_n|)$, is the
``norm'' of the Lorentz space $\ell^\pinf$.  The cardinality of a set
$A$ is denoted by $|A|$.

\section{Definition of the space $E$ and simple properties}

Let $D = \{(i,j): i,j \in \N\, , i \geq j\}$ and let $G$ be the vector
lattice of all functions $x:D\to \R\ $ having finite support.
Then let
\[B =
   \{b=(b_i) \in c_{00}: ib_i \in \N\cup \{0\} \mbox{ for all } i,
    \|b_i\|_2 \leq 1 \}.          \]
For all $b \in B$, define $x_b \in G$ by
\[ x_b(i,j) = \left\{ \begin{array}{ll}
                       1 & 1 \leq j \leq ib_i, i \in \N \\
                       0 & \mbox{otherwise}.
                      \end{array}
              \right. \]
Let $U$ be the convex solid hull of $\{x_b: b \in B\}$.  Define a
seminorm $\rho$ on $G$ by
\[ \rho(x) =
\lnorm\biggl(\frac{1}{i}
\sum^{i}_{j=1}|x(i,j)|\biggr)^\infty_{i=1}\rnorm_2
.
\]
Two elements $x, y \in G$ are {\em row disjoint} if
\[ \sumi |x(i,j)| \cdot \sumi |y(i,j)| = 0 \]
for all $i$.  Fix $1 < p < 2$, and let
\begin{eqnarray*}
 A & = & \{y_1 + \cdots + y_m: m \in \N\, ,
y_1, \ldots , y_m \in
U\ \mbox{pairwise row disjoint}, \\
& & \|(\rho(y_1),\ldots,\rho(y_m))\|_{p,\infty} \leq 1 \} .
\end{eqnarray*}
For $y = y_1+\cdots+y_m$, where $y_1,\ldots,y_m$ are as above, we say
that the sum on the right is a {\em representative} of the element $y
\in A$, and $m$ is the {\em length} of the representative.
Furthermore, let
\[
 \phi(y) = \min\{\|(\rho(y_1),\ldots,\rho(y_m))\|_{\infty}:
     y_1+\ldots+y_m  \mbox{ is a representative of } y\}
\]
for all $y \in A$. The minimum exists since $y$ is finitely supported.
An element $y \in A$ with $\phi(y) \leq \epsilon$ will be called
$\epsilon${\em -small}.  A subset of $A$ is $\epsilon$-{\em small} if
all of its members are $\epsilon$-small.
Finally, let $V$  be the convex hull of $A$. \\

It is easy to see that $V$ is a convex solid subset of $G$ such that
$\cap_{\lambda>0}\lambda V = \{0\}$. Hence the gauge functional $\tau$
of $V$ is a lattice norm on $G$.  Let $E$ be the completion of $G$
with respect to the norm $\tau$.  It is clear that $E$ has an
unconditional basis.

\begin{lem}
Let $C = \sqrt{\sum^\infty_{n=1}n^{-2/p}}$.  Then $\rho(x) \leq C$
for all $x
\in V$.
\end{lem}

\begin{pf}
Let $y = \sum^m_{i=1}y_i \in A$,
where $y_1, \ldots , y_m \in U$ are pairwise row disjoint and
\[ \|(\rho(y_1),\ldots,\rho(y_m))\|_{p,\infty} \leq 1 . \]
Then
\[ \|(\rho(y_1),\ldots,\rho(y_m))\|_2 \leq C . \]
Since  $y_1, \ldots , y_m \in U$ are pairwise row disjoint,
\[ \rho(y) =  \|(\rho(y_1),\ldots,\rho(y_m))\|_2 \leq C . \]
The result now follows easily since $V = co(A)$.
\end{pf}

\begin{pr}\label{quot}
$E$ has a quotient space isomorphic to $\ell^2$.
\end{pr}

\begin{pf}
This is equivalent to showing that $E' = G'$ contains an isomorph of
$\ell^2$.  For all $i \in \N$, define $z_i:G \to \R$ by
$\langle x,z_i\rangle = i^{-1}\sumi x(i,j).$  We will show that the
sequence $(z_i)$ in $G'$ is equivalent to the $\ell^2$-basis.   Let
$(b_i)$ be a
finitely supported
sequence on the unit
sphere of $\ell^2$. For any $x \in V$,
\begin{eqnarray*}
\langle x, \sum b_iz_i \rangle & = & \sum_i\frac{b_i}{i}\sumi x(i,j)
\\
 & \leq &
  \lnorm\biggl(\frac{1}{i}
\sum^{i}_{j=1}x(i,j)\biggr)^\infty_{i=1}\rnorm_2 \\
 & \leq & \rho(x) \\
 & \leq & C .
\end{eqnarray*}
Hence $\tau{'}(\sum b_iz_i) \leq C$, where $\tau{'}$
denotes the norm dual
to $\tau$.\\
\indent On the other hand, we claim that  $\tau{'}(\sum b_iz_i) \geq
2/9$.  Indeed, if $\|(b_1, \ldots, b_{4})\|_2 \geq
\sqrt{5}/3$, then since $(z_i)$ is clearly pairwise disjoint (in the
lattice sense) and $\tau'(z_i) \geq 1$ for all $i$,
\begin{eqnarray*}
 \tau{'}(\sum b_iz_i) & \geq & \sup |b_i| \\
 & \geq & \frac{1}{2}\|(b_1,\ldots,b_{4})\|_2 \\
 & \geq & \sqrt{5}/6 \\
 & \geq & 2/9.
\end{eqnarray*}
Now if  $\|(b_1, \ldots, b_{4})\|_2 <
\sqrt{5}/3$, then $\alpha \equiv
\|(b_{5}, b_{6}, \ldots)\|_2 >
2/3$.  For all $i > 4$, let $m_i$ be the
largest non-negative integer
$\leq i|b_i|$.  Then $m \equiv (0,0,0,0,m_{5}/5,
m_{6}/6,\ldots) \in B$.  Let $y \in G$ be given by
\[ y(i,j) = \left\{ \begin{array}{ll}
            \sgn b_i & 1 \leq j \leq m_i\ , i > 4 \\
                       0 & \mbox{otherwise}.
                      \end{array}
              \right. \]
Then $|y| = x_m \in U$, and hence $y \in V$.  Therefore,
\begin{eqnarray*}
\tau'(\sum b_iz_i) & \geq &
\langle y, \sum b_iz_i \rangle \\
 & = &   \sum_i\frac{b_i}{i}\sumi y(i,j) \\
 & = & \sumsv\frac{m_i|b_i|}{i} . \\
\end{eqnarray*}
Since $m_i/i \geq |b_i| - 1/i$ for $ i > 4$, we have
\begin{eqnarray*}
\tau'(\sum b_iz_i) & \geq & \sumsv |b_i|^2
- \sumsv \frac{|b_i|}{i} \\
 & \geq & \alpha^2 - \alpha\sumsv i^{-2} \\
 & \geq & \alpha^2 - \frac{\alpha}{3} \\
 & > & \frac{2}{9}\, ,
\end{eqnarray*}
since $\alpha > 2/3$.
\end{pf}

Let us try to provide some motivation behind the
particular construction of the set $V$.  First off, we want to fix
things so that the sequence $(z_i)$ as given in
Proposition \ref{quot}
is equivalent to the $\ell^2$-basis. This has to be
done without
introducing any sequence biorthogonal to $(z_i)$ whose linear
combinations can be ``normed'' by vectors in $\spn\{z_i\}$.
Thus we
settled on using the set $\{x_b: b\in B\}$ to ``norm'' the
vectors in
$\spn\{z_i\}$.  This puts all the elements of the set
$U$ into the unit
ball of the space $E$ we are trying to construct.
(The word ``solid'' is
needed since we want a space with an unconditional basis.)
To try to
obtain
$c_0$-saturation, we admit into the unit ball of $E$ {\em some}
elements of the
form $y_1+\cdots+ y_m$, where $y_1,\ldots, y_m \in U$ are
pairwise row
disjoint. However, to keep the equivalence of $(z_i)$ to the
$\ell^2$-basis,
the elements admitted must have uniformly bounded $\rho$-norms.
The $\ell^{p,\infty}$-norm (which appears in the definition
of the set
$A$) is chosen because it majorizes the $\ell^2$-norm,
and the closed
subspace of $\ell^{p,\infty}$ generated by the unit
vector basis is
$c_0$-saturated.

\section{Proof that $E$ is $c_0$-saturated}

Let $k \in \N$, a collection of real sequences is $k${\em -disjoint}
if the pointwise product of any $k+1$ members of the collection is the
zero sequence; equivalently, if at most $k$ of the sequences can be
non-zero at any fixed coordinate.  We begin with some elementary
lemmas.

\begin{lem}\label{elem}
Let $\{y_1,\ldots, y_k\}$ be a finite subset of $\ell^2$, then
\[ \lnorm\sum^k_{i=1}y_i\rnorm_2 \leq
\sqrt{k}\biggl(\sum^k_{i=1}\|y_i\|^2_2\biggr)^{1/2} . \]
\end{lem}

\begin{lem}\label{kdisj}
For any $k,n \in$ {\rm \N}, and any $k$-disjoint subset $\{x_1,
\ldots, x_n\}$
of the unit ball of $\ell^2$,
\[ \lnorm\sum^n_{i=1} x_i\rnorm_2 \leq \sqrt{kn} . \]
\end{lem}

\begin{pf}
Write $x_i = (x_i(j))$ for $1 \leq i \leq n$.  For each $j$, let
$y_1(j),\ldots, y_n(j)$ be the decreasing rearrangement of
$|x_1(j)|,\ldots, |x_n(j)|$.  Then let $y_i = (y_i(j))$.  By the
$k$-disjointness, $y_i = 0$ for $i > k$.  Hence $|\sum^n_{i=1}x_i|
\leq \sum^k_{i=1}y_i$ (pointwise order).  Therefore,
\begin{eqnarray*}
\lnorm\sum^n_{i=1}x_i\rnorm_2 & \leq &
\lnorm\sum^k_{i=1}y_i\rnorm_2 \\
 & \leq & \sqrt{k}\biggl(\sum^k_{i=1}\|y_i\|^2_2\biggr)^{1/2}
\hspace{1em}
\mbox{by Lemma \ref{elem}} \\
 & = & \sqrt{k}\biggl(\sum^k_{i=1}\sum_jy_i(j)^2\biggr)^{1/2} \\
 & = & \sqrt{k}\biggl(\sum^n_{i=1}\sum_jx_i(j)^2\biggr)^{1/2} \\
 & = & \sqrt{k}\biggl(\sum^n_{i=1}\|x_i\|^2_2\biggr)^{1/2} \\
 & \leq & \sqrt{kn} .
\end{eqnarray*}
\end{pf}

\begin{lem}\label{small}
If $x \in U$ satisfies $\|x\|_\infty \leq
\epsilon$ for some  $\epsilon > 0$, then $\rho(x) \leq
\sqrt{\epsilon}$.
\end{lem}

\begin{pf}
First of all, for any $\delta > \epsilon$, one can find $z \in
co\{x_b: b \in B\}$ such that $|x| \leq z$ and $\|z\|_\infty <
\delta$. By a small perturbation, it may even be assumed that $z$ is
a convex combination with rational coefficients.  Then, $z$ can be
expressed in the form
$z = n^{-1}\sum^n_{i=1}x_{b_i}$, where $b_i \in B, 1
\leq i \leq n.$  Let $j$ be the greatest integer $\leq \delta n$.
Then $\|z\|_\infty \leq \delta$ implies $\{b_1,\ldots, b_n\}$ is
$j$-disjoint.  Therefore
\[ n^{-1}\lnorm\sum^n_{i=1}b_i\rnorm_2 \leq \sqrt{\frac{j}{n}} \leq
\sqrt{\delta}  \]
by Lemma \ref{kdisj}.  It remains to observe that the leftmost
quantity in the above equation is precisely $\rho(z)$, and $\delta >
\epsilon$ is arbitrary.
\end{pf}

The next lemma is a quantitative version of the fact that the unit
vector basis of $\ell^\pinf$ generates a $c_0$-saturated closed
subspace.

\begin{lem}\label{elpinf}
Let $(a_i)$ be a real sequence, and $0 = n_0 < n_1 < \ldots$ a
sequence of integers so that
\begin{enumerate}
\item $\|(a_{n_k+1},\ldots,a_{n_{k+1}})\|_\pinf
         \leq 1$, and
\item $\|(a_{n_k+1},\ldots,a_{n_{k+1}})\|_\infty \leq
         n_k^{-1/p}$
\end{enumerate}
for all $k \geq 0$.  Then $\|(a_1, a_2, \ldots)\|_\pinf \leq 2$.
\end{lem}

\begin{pf}
Assume the contrary. Then there exists $n$ such that $a^*_n >
2n^{-1/p}$. Hence $J = \{i: |a_i| > 2n^{-1/p} \}$ has cardinality
$\geq n$.  Since $a_i \to 0$ as $i \to \infty$, $J$ is finite.  Let
$j$ be the largest element in $J$.  Choose $k$ so that $n_k < j \leq
n_{k+1}$. Then
\[ 2n^{-1/p} < |a_j| \leq n_k^{-1/p} . \]
Hence $n_k < 2^{-p}n$.  Therefore,
\[ |J \cap \{n_k+1, \ldots, n_{k+1}\}| > n(1-2^{-p}). \]
Consequently,
\begin{eqnarray*}
\|(a_{n_k+1}, \ldots, a_{n_{k+1}})\|_\pinf & \geq &
  2n^{-1/p}\|(\overbrace{1, \ldots, 1}^{n(1-2^{-p})})\|_\pinf \\
 & > & 1 ,
\end{eqnarray*}
a contradiction.
\end{pf}

To prove that $E$ is $c_0$-saturated, it suffices to show that every
pairwise row disjoint sequence $(y_i) \subseteq V$ which satisfies
$\|y_i\|_\infty \to 0$ has a subsequence $(y_{i_k})$ with the
property  $\sup_n\tau(\sum^n_{k=1}y_{i_k}) < \infty$.  To begin with,
consider  the simple case when $(y_i) \subseteq A$.  Pick a
representative $y_i = y_i(1)+\ldots+y_i(m_i)$ for each $i$.  Clearly,
$\sup_{1\leq j\leq m_i}\!\|y_i(j)\|_\infty \to 0$ as $i \to \infty$.
Hence, by Lemma \ref{small}, $\sup_{1\leq j\leq m_i}\!\rho(y_i(j))
\to 0$ as well.  By dropping to a subsequence, we may assume that
\[ \|(\rho(y_i(1)),\ldots,\rho(y_i(m_i)))\|_\infty \leq
         \biggl(\sum^{i-1}_{k=1}m_k\biggr)^{-1/p}. \]
Then, by Lemma \ref{elpinf}, for all $i$, $2^{-1}\!\sum^i_{k=1}y_k$
has a
representative $2^{-1}\!\sum^i_{k=1}\!\sum^{m_k}_{j=1}y_k(j) \in A$.
Therefore, $\tau(\sum^i_{k=1}y_k) \leq 2$ for all $i$, as required.
To tackle the general case where we know only that $(y_i) \subseteq
V$, we introduce the following definition.\\

\noindent{\bf Definition.} A sequence $(y_i) \subseteq V$ is
called {\em
strongly decreasing} if there exists
a sequence $(\epsilon_i)$ decreasing to $0$ such that
for every $i$, there is a $\epsilon_i$-small subset $A_i$
of $A$ with $y_i \in co(A_i)$.  \\

Basically the same argument as above shows that

\begin{pr}\label{strongdec}
A $\tau$-normalized, pairwise row disjoint, strongly decreasing
sequence in $V$ has a subsequence equivalent to the $c_0$-basis.
\end{pr}

What will be shown is that every sequence $(y_i)$ in $V$ with
$\|y_i\|_\infty
\to 0$ can be written as the sum of a strongly decreasing sequence
and
a $\tau$-null sequence.  The
proof of this fact relies on the

\begin{mainlem}
If $x \in U$ satisfies $\|x\|_\infty \leq \epsilon$ for some
$\epsilon
> 0$, then $\tau(x) \leq 5\epsilon^{1/4}$.
\end{mainlem}

The proof of the Main Lemma will be given in the next section.
Assuming the result, we continue with

\begin{lem}
Let $y \in V$ satisfy $\|y\|_\infty \leq \epsilon$, then there exist
 a
 $\epsilon^{1/8p}$-small subset $S$ of $A$, and  $u \in co(S)$,
such that $\tau(y-u) \leq 5\epsilon^{1/8}$.
\end{lem}

\begin{pf}
There is no loss of generality in assuming that $y \geq 0$.
Express $y$ as a convex combination $\sum^j_{i=1}\alpha_iy_i$ where
$y_i \in A,\, 1 \leq i \leq j$.  Since $A$ is solid, it may be
assumed that $y_i \geq 0$ for all $i$.  Choose a
representative $\sum^m_{l=1}y_i(l)$ for each $y_i$. (All the
representatives can be made to have the same length $m$ by adding on
zeros if necessary.)  Note that $y_i(l) \geq 0$ for all $i$ and $l$.
Let $\delta = \epsilon^{1/8p}$ and let $A_i = \{l: \rho(y_i(l)) >
\delta\}$.  Since $y \in V$,
\begin{eqnarray*}
1 & \geq & \|(\rho(y_i(1)), \ldots, \rho(y_i(m)))\|_{p,\infty} \\
 & \geq & \delta\|(\underbrace{1, \ldots, 1}_{|A_i|})\|_{p,\infty} \\
 & = & \delta|A_i|^{1/p} .
\end{eqnarray*}
Hence $|A_i| \leq \delta^{-p}$ for all $i$.  Let $r$ be the greatest
integer $\leq \delta^{-p}$.  Relabeling, we may assume that each $A_i$
is an initial interval $\{1, \ldots, r_i\}$, where $r_i \leq r$.
Define $v_l = \sum^j_{i=1}\alpha_iy_i(l), 1 \leq l \leq r$, and let $v
= \sum^r_{l=1}v_l$.  Then $v_l \in U,$ and $\|v_l\|_\infty \leq
\|y\|_\infty \leq \epsilon$ for all $l$.  By the Main Lemma,
we obtain
\begin{eqnarray*}
\tau(v) & \leq & \sum^r_{l=1}\tau(v_l) \\
 & \leq & 5\epsilon^{1/4}r \\
 & \leq & 5\epsilon^{1/4}\delta^{-p} \\
 & = & 5\epsilon^{1/8} .
\end{eqnarray*}
Now let $u_i = \sum^m_{l=r+1}y_i(l)$, $S = \{u_i: 1 \leq i \leq j\},$
and $u = \sum^j_{i=1}\alpha_iu_i$. It is clear that $S$ is a
$\epsilon^{1/8p}$-small subset of $A$, $u \in co(S)$, and $\tau(y-u) =
\tau(v) \leq 5\epsilon^{1/8}.$
\end{pf}

The promised result is now immediate.

\begin{pr}\label{split}
Every sequence $(y_i)$ in $V$ satisfying $\|y_i\|_\infty \to 0$ can be
written as the sum of a strongly decreasing sequence and a $\tau$-null
sequence.
\end{pr}

Using Propositions \ref{strongdec}, \ref{split}, and standard
arguments, we obtain

\begin{thm}
Every normalized, pairwise disjoint, finitely supported sequence
$(y_i)$ in $E$ with $\|y_i\|_\infty \to 0$ has a subsequence
equivalent to the $c_0$-basis.  Consequently, $E$ is $c_0$-saturated.
\end{thm}

\section{Proof of the Main Lemma}

In this section, we prove the Main Lemma.  In fact, we will show that
if $x \in U$ satisfies $\|x\|_\pinf \leq \epsilon$,
then $x \in 5\ep\! A$.
The basic idea is as follows.  Given such a $x$,
Lemma \ref{small} says
that $\rho(x) \leq \sqrt{\epsilon}$.  Let $n$ be any natural number
$\leq 5^p\epsilon^{-p/4}$.  Then if $x$ is written as a pairwise row
disjoint sum $x = \sum^n_{i=1}x_i$,
\begin{eqnarray*}
\|(\rho(x_1), \ldots, \rho(x_n))\|_\pinf & \leq &
  \rho(x)\|(\overbrace{1, \ldots, 1}^n)\|_\pinf \\
 & \leq & \sqrt{\epsilon}n^{1/p} \\
 & \leq &  5\ep .
\end{eqnarray*}
What we need to show is that the $x_i$'s can be chosen so that they
come from a small multiple ($5\ep$) of $U$.  The proof of the next
Lemma is left to the reader.

\begin{lem}\label{select}
Let $(a_i)^l_{i=1}$ be numbers in $[0,1]$ such that $\sum^l_{i=1}a_i
> 1$, then there exists $S \subseteq \{1, \ldots, l\}$ such that
$1/2 \leq \sum_{i\in S}a_i \leq 1$.
\end{lem}

For a real-valued matrix $a = (a_{i,j})^{k\hspace{1em}l}_{i=1 j=1}$,
we
define $\Sigma(a) = \sum_{i,j}a_{i,j}$, and, for each $j$,
$s_j(a) = \min\{i: a_{i,j} \neq 0\}\ (\min\emptyset = 0)$.

\begin{lem}\label{partition}
Let $a = (a_{i,j})^{k\hspace{1em}l}_{i=1 j=1}$ be a $[0,1]$-valued
matrix
such that $\Sigma(a) \leq  M$.
Then there is a partition $R_1, \ldots, R_n$ of $R = \{(i,j): 1 \leq i
\leq k, 1 \leq j \leq l\}$ such that
\begin{enumerate}
\item $n \leq 2M + k$,
\item $\displaystyle{\sum_{(i,j)\in R_m}\! a_{i,j}} \leq 1,
  \hspace{2em} 1 \leq m \leq n, $ and
\item $|R_m \cap \{(i,j): 1 \leq i \leq k\}| \leq 1$ for all $ 1 \leq
m \leq n$ and $1 \leq j \leq l$.
\end{enumerate}
\end{lem}

\begin{pf}
A $[0,1]$-valued matrix $b = (b_{i,j})^{k\hspace{1em}l}_{i=1 j=1}$
will be called {\em reducible} if $\sum^l_{j=1}\!b_{s_j(b),j} > 1$,
where $b_{s_j(b),j}$ is taken to be $0$ if $s_j(b) = 0$.  For a
reducible matrix $b$, Lemma \ref{select} provides a set $S(b)
\subseteq \{(s_j(b),j): s_j(b) > 0\}$ such that
$1/2 \leq \sum_{(i,j)\in S(b)}\! b_{i,j} \leq 1$.  We also let the
{\em reduced} matrix $b'$ be given by
\[ b'_{i,j} = \left\{ \begin{array}{ll}
                       b_{i,j} & \mbox{if $(i,j) \notin S(b)$} \\
                       0       & \mbox{otherwise}
                      \end{array}
              \right. \]
Now let the matrix $a$ be as given.
Since all the conditions are invariant with respect to permutations of
the entries of $a$ within columns, we may assume that $a_{i,j} \geq
a_{i+1,j}$. Let $a^1 = a$.  If $a^1$ is reducible, let $R_1 = S(a^1)$
and $a^2 = {a^1}'$.  Inductively, if $a^r$ is reducible, let $R_r =
S(a^r)$ and $a^{r+1} = {a^r}'$.  Note that $\Sigma(a^{r+1}) \leq
\Sigma(a^r) - 1/2$ if $a^r$ is reducible.  Therefore, there must be a
$t < 2M$ such that $a^{t+1}$ is not reducible.  Now let
\[ R_{t+i} = \{(i,j): (i,j) \notin R_1 \cup \ldots \cup R_t \} \]
for $1 \leq i \leq k$.  Let $n = t + k$.  It is clear that the
collection $R_1, \ldots, R_n$ satisfies the requirements.
\end{pf}

\noindent{\em Proof of the Main Lemma:} Since $x \in U$ already
implies
$\tau(x) \leq 1$, we may assume $\epsilon < 1$.
As in the proof of Lemma \ref{small},
it suffices to prove the Main Lemma for those $x$'s which have the
form $x = M^{-1}\!\sum^M_{i=1}x_{b_i}$. Write $b_i = (b_{i,j})_j$ for
$1 \leq i \leq M$.  Since $b_i \in c_{00}$, there exists $l$ such that
$b_{i,j} = 0$ for all $j > l$.  Also, $\|x\|_\infty \leq \epsilon$
implies
$\{b_1, \ldots, b_M\}$ is $k$-disjoint, where $k$ is the greatest
integer $\leq \epsilon M$.
For each $j$, let $a_{1,j}, \ldots,
a_{M,j}$ be the decreasing rearrangement of $b^2_{1,j}, \ldots,
b^2_{M,j}$.  By the $k$-disjointness, $a_{i,j} =  0$ for all $i > k$.
Now let $a = (a_{i,j})^{k\hspace{1em}l}_{i=1 j=1}$.  Then $a$ is a
$[0,1]$-valued matrix, and
\[ \Sigma(a) = \sum^n_{i=1}\!\sum^{l}_{j=1}b^2_{i,j}
  = \sum^n_{i=1}\|b_i\|^2_2 \leq M .\]
Let $\eta = (2\epsilon^{-p/4}-1)^{-1}$.  Since we are assuming that
$\epsilon < 1$, $\eta$ is positive and $\leq \ep$.
Choose the smallest integer $j_1$ such
that $\sumk\!\sum^{j_1}_{j=1}\!a_{i,j} > \eta M$.
Since $\sumk\!a_{i,j}
\leq k \leq \epsilon M$ for all $j$, $\sumk\!\sum^{j_1}_{j=1}\!a_{i,j}
\leq (\eta + \epsilon)M$.  If $\sumk\sum^{l}_{j=j_1+1}a_{i,j} > \eta
M$, we choose the smallest integer $j_2$ such that
$\sumk\!\sum^{j_2}_{j=j_1+1}\!a_{i,j} > \eta M$.
Continuing inductively,
we obtain $0 = j_0 < j_1 < \ldots < j_t = l$ such that
\[ \eta M < \sumk\sum^{j_{m+1}}_{j=j_m+1}\!a_{i,j} \leq (\eta
+\epsilon)M \]
for $0 \leq m \leq t-2$, and
\[ \sumk\sum^{j_t}_{j=j_{t-1}+1}\!a_{i,j} \leq \eta M .\]
Note that
\[ (t-1)\eta M < \sum^{t-2}_{m=0}\sumk\sum^{j_{m+1}}_{j=j_m+1}a_{i,j}
  \leq \Sigma(a) \leq M \]
implies $t < \eta^{-1}+1$.  With $M$ replaced by $(\eta + \epsilon)M$,
apply Lemma \ref{partition} to each submatrix
$(a_{i,j})^{k\hspace{1em}j_{m+1}}_{i=1 j=j_m+1}$
to obtain a partition
$R^m_1, \ldots, R^m_{r_m}$ of
$\{(i,j): 1 \leq i \leq k,~j_m~<~j~\leq~j_{m+1}\}\\
(0 \leq m < t)$ such that
\begin{itemize}
\item $r_m \leq 2(\eta + \epsilon)M + k \leq (2\eta + 3\epsilon)M$,
\item $\displaystyle{\sum_{(i,j)\in R^m_\nu}\!a_{i,j}} \leq 1$,
\hspace{1em} $1 \leq \nu \leq r_m$, and
\item $|R^m_\nu \cap \{(i,j): 1 \leq i \leq k\}| \leq 1$ for all $m$
and $\nu$.
\end{itemize}
Finally, if $j$ is such that $|R^m_\nu \cap \{(i,j): 1 \leq i \leq
k\}| = 1$, and $(i,j)$ is the unique element of this set, let
$d^m_\nu(j) = \sqrt{a_{i,j}}$; otherwise, let $d^m_\nu(j) = 0$.  Then
for all $m$ and $\nu$, the sequence $d^m_\nu =
(d^m_\nu(j))^\infty_{j=1} \in B$.  The following points are worth
noting.
\begin{itemize}
\item For every $j$, the nonzero numbers in the list
$(d^m_\nu(j))^{r_m\hspace{.5em}t-1}_{\nu=1 m=0}$ is a rearrangement of
the nonzero numbers in the list $(b_{1,j}, \ldots, b_{M,j})$,
\item If $c_i = (c_{i,j})_j \in B, 1 \leq i \leq \alpha$
and $(i_0,j_0) \in D$, then
\[ (x_{c_1} + \ldots + x_{c_\alpha})(i_0,j_0) =
  |\{i: i_0c_{i,i_0} \geq j_0\}| . \]
\end{itemize}
Using these, it is easy to see that
\[ \sum^M_{i=1}x_{b_i} =
\sum^{t-1}_{m=0}\sum^{r_m}_{\nu=1}x_{d^m_\nu}. \]
Now let $y_m =
M^{-1}\sum^{r_m}_{\nu=1}x_{d^m_\nu}$, $0 \leq m < t$.  Then $y_0,
\ldots, y_{t-1}$ are pairwise row disjoint, and $x = y_0 + \cdots +
y_{t-1}$. Also,
\begin{eqnarray*}
y_m & = & \frac{r_m}{M}(r^{-1}_m\sum^{r_m}_{\nu=1}x_{d^m_\nu}) \\
   & \in & \frac{r_m}{M}U \\
   & \subseteq & (2\eta + 3\epsilon)U \\
   & \subseteq & 5\ep U.
\end{eqnarray*}
Furthermore,
\begin{eqnarray*}
\|(\rho(y_0), \ldots, \rho(y_{t-1}))\|_\pinf
 & \leq & \rho(x)\|(\overbrace{1,\ldots,1}^t)\|_\pinf \\
 & \leq & \epsilon^{1/2}t^{1/p}  \\
 & \leq & 2\ep .
\end{eqnarray*}
Therefore, $x \in 5\ep A$, so $\tau(x) \leq 5\ep$,
as required. $\Box$\\

\noindent{\bf Remark\ }  Let $(e_i)$ be an unconditional basis of a
Banach space $F$.  If $(e_i)$ has the property
\[ (*) \left\{
\begin{array}{l}
\mbox{every normalized block basis $(\sum^{j_{k+1}}_{i=j_k+1}a_ie_i)$
of $(e_i)$ such that} \\
\mbox{$a_i \to 0$ has a subsequence equivalent to the
$c_0$-basis,}
\end{array}
\right. \]
then $F$ is $c_0$-saturated.  Since the unit vector basis of the space
$E$ has property ($*$), we see that it alone is not enough to insure
that $F'$ is $\ell^1$-saturated.  However, the following proposition
is easy to obtain.

\begin{pr}
Let $(e_i)$ be an unconditional basis of a Banach space $F$.  If
$(e_i)$ has property {\rm ($*$)}, and $\|\sum_{i\in A_k}e_i\| \to
\infty$
whenever $(A_k)$ is a sequence of subsets
of\/ {\rm \N\,} such that $\max A_k < \min A_{k+1}$ and $|A_k|
\to \infty$,
then $F'$ is $\ell^1$-saturated.
\end{pr}

\flushleft
\vspace{.5in}
Department of Mathematics\\National University of Singapore\\
Singapore 0511\\ e-mail(bitnet) : matlhh@nusvm

\end{document}